# OPTIMAL PREDICTIVE MODEL SELECTION[1]

By Maria Maddalena Barbieri and James O. Berger

*Università Roma Tre and Duke University*

Often the goal of model selection is to choose a model for future prediction, and it is natural to measure the accuracy of a future prediction by squared error loss. Under the Bayesian approach, it is commonly perceived that the optimal predictive model is the model with highest posterior probability, but this is not necessarily the case. In this paper we show that, for selection among normal linear models, the optimal predictive model is often the *median probability model*, which is defined as the model consisting of those variables which have overall posterior probability greater than or equal to 1/2 of being in a model. The median probability model often differs from the highest probability model.

**1. Introduction.** Consider the usual normal linear model

$$\mathbf{y} = \mathbf{X}\boldsymbol{\beta} + \boldsymbol{\varepsilon}, \tag{1}$$

where $\mathbf{y}$ is the $n \times 1$ vector of observed values of the response variable, $\mathbf{X}$ is the $n \times k$ $(k < n)$ full rank design matrix of covariates, and $\boldsymbol{\beta}$ is a $k \times 1$ vector of unknown coefficients. We assume that the coordinates of the random error vector $\boldsymbol{\varepsilon}$ are independent, each with a normal distribution with mean 0 and common variance $\sigma^2$ that can be known or unknown. The least squares estimator for this model is thus $\widehat{\boldsymbol{\beta}} = (\mathbf{X}'\mathbf{X})^{-1}\mathbf{X}'\mathbf{y}$.

Equation (1) will be called the *full* model, and we consider selection from among submodels of the form

$$M_{\mathbf{l}} : \mathbf{y} = \mathbf{X}_{\mathbf{l}}\boldsymbol{\beta}_{\mathbf{l}} + \boldsymbol{\varepsilon}, \tag{2}$$

where $\mathbf{l} = (l_1, l_2, \ldots, l_k)$ is the model index, $l_i$ being either 1 or 0 as covariate $x_i$ is in or out of the model (or, equivalently, if $\beta_i$ is set equal to zero); $\mathbf{X}_{\mathbf{l}}$

Received January 2002; revised April 2003.

[1]Supported by NSF Grants DMS-98-02261 and DMS-01-03265, and by Consiglio Nazionale delle Ricerche and Ministero dell'Istruzione, dell'Università e della Ricerca.

contains the columns of $\mathbf{X}$ corresponding to the nonzero coordinates of $\mathbf{l}$; and $\boldsymbol{\beta}_{\mathbf{l}}$ is the corresponding vector of regression coefficients.

Upon selecting a model, it will be used to predict a future observation

$$(3) \qquad y^* = \mathbf{x}^* \boldsymbol{\beta} + \varepsilon,$$

where $\mathbf{x}^* = (x_1^*, \ldots, x_k^*)$ is the vector of covariates at which the prediction is to be performed. The loss in predicting $y^*$ by $\widehat{y}^*$ will be assumed to be the squared error loss

$$(4) \qquad L(\widehat{y}^*, y^*) = (\widehat{y}^* - y^*)^2.$$

With the Bayesian approach to model selection it is commonly perceived that the best model will be that with the highest posterior probability. This is true under very general conditions if only two models are being entertained [see Berger (1997)] and is often true in the variable selection problem for linear models having orthogonal design matrices [cf. Clyde (1999) and Clyde and George (1999, 2000)], but is not generally true. Indeed, even when only three models are being entertained essentially nothing can be said about which model is best if one knows only the posterior probabilities of the models. This is demonstrated in Section 5, based on a geometric representation of the problem.

For prediction of a single $y^*$ at a specific $\mathbf{x}^*$, one can, of course, simply compute the posterior expected predictive loss corresponding to each model and choose the model that minimizes this expected loss. In such a scenario, however, choosing a specific model makes little sense; one should, rather, base the prediction on Bayesian *model averaging* [cf. Clyde (1999) and Hoeting, Madigan, Raftery and Voliksky (1999)]. The basic use of *model selection* for prediction is when, because of outside constraints, a single model must be selected for repeated use in future predictions. (Note that we are assuming that these constraints preclude use of the Bayesian model averaging answer.) It is natural to assume that these future predictions will be made for covariates $\mathbf{x}^*$ that arise according to some distribution. We further assume that the $k \times k$ expectation matrix corresponding to this distribution,

$$(5) \qquad \mathbf{Q} = \mathbf{E}[(\mathbf{x}^*)'(\mathbf{x}^*)],$$

exists and is positive definite. A frequent choice is $\mathbf{Q} = \mathbf{X}'\mathbf{X}$, which is equivalent to assuming that the covariates that will occur in the future are like those that occurred in the data. (Strictly, this would yield $\mathbf{Q} = \frac{1}{n}\mathbf{X}'\mathbf{X}$, but constants that are common across models can be ignored.)

In this scenario, one could still simply compute the expected predictive loss corresponding to each model and minimize, but the expectation would now also be over $\mathbf{x}^*$. This can add quite a computational burden, especially when there are many models to consider. Bayesian MCMC schemes have been developed that can effectively determine the posterior model probabilities $P(M_{\mathbf{l}}|\mathbf{y})$, but adding an expectation over $\mathbf{x}^*$ and a minimization over $\mathbf{l}$ can be prohibitive [although see Müller (1999)]. We thus sought to determine if there are situations in which it is possible to give the optimal predictive model solely in terms of the posterior model probabilities.

Rather general characterizations of the optimal model turn out to be frequently available but, quite surprisingly, the characterizations are not in terms of the highest posterior probability model, but rather in terms of what we call the *median probability model*.

DEFINITION 1. The *posterior inclusion probability* for variable $i$ is

$$p_i \equiv \sum_{\mathbf{l}:\, l_i = 1} P(M_{\mathbf{l}}|\mathbf{y}), \tag{6}$$

that is, the overall posterior probability that variable $i$ is in the model.

If it exists, the *median probability model* $M_{\mathbf{l}^*}^*$ is defined to be the model consisting of those variables whose posterior inclusion probability is at least $1/2$. Formally, $\mathbf{l}^*$ is defined coordinatewise by

$$l_i^* = \begin{cases} 1, & \text{if } p_i \geq \frac{1}{2}, \\ 0, & \text{otherwise.} \end{cases} \tag{7}$$

It may happen that the set of covariates defined by (7) does not correspond to a model under consideration, in which case the median probability model will not exist. There are, however, two important cases in which the median probability model is assured to exist. The first is in the problem of variable selection when any variable can be included or excluded from the model (so that all vectors $\mathbf{l}$ are possible).

The second case of particular interest is when the class of models under consideration has a graphical model structure.

DEFINITION 2. Suppose that for each variable index $i$ there is a corresponding index set $I(i)$ of other variables. A subclass of linear models has *graphical model structure* if it consists of all models satisfying the condition "for each $i$, if variable $x_i$ is in the model, then variables $x_j$ with $j \in I(i)$ are in the model."

It is straightforward to show that if a subclass of linear models has graphical model structure, then the median probability model will satisfy the condition in the definition and, hence, will itself be a member of the subclass.

One common example of a class of linear models having graphical model structure is the class of all models that can be constructed from certain main effects and their interactions up to a certain order, subject to the condition that if a high order interaction of variables is in a model then all lower order interactions (and main effects) of the variables must be in the model.

A second example of a subclass having graphical model structure is a sequence of nested models,

$$(8) \qquad M_{\mathbf{l}(j)}, \qquad j = 0, \dots, k, \text{ where } \mathbf{l}(j) = (1, \dots, 1, 0, \dots, 0),$$

with $j$ ones and $k - j$ zeroes. Examples of this scenario include polynomial regression, in which $j$ refers to the polynomial order used, and autoregressive time series, in which $j$ refers to the allowed lag. Note that for nested models the median probability model has a simpler representation as $M_{\mathbf{l}(j^*)}$, where $j^*$ is such that

$$(9) \qquad \sum_{i=0}^{j^*-1} P(M_{\mathbf{l}(i)}|\mathbf{y}) < \tfrac{1}{2} \quad \text{and} \quad \sum_{i=0}^{j^*} P(M_{\mathbf{l}(i)}|\mathbf{y}) \geq \tfrac{1}{2}.$$

In other words, one just lists the sequence of posterior model probabilities and sums them up until the sum exceeds 1/2. The model at which the exceedance occurs is the median probability model.

The above special cases also define the scenarios that will be investigated in this paper. The goal will be to provide conditions under which the median probability model is the optimal predictive model. The conditions are primarily restrictions on the form of the predictors for $y^*$. The restrictions are fairly severe, so that the results can best be thought of as applying primarily to default Bayes or empirical Bayes types of procedures.

Initially we had sought to find conditions under which the highest posterior probability model was the optimal predictive model. It came as quite a surprise to find that any optimality theorems we could obtain were, instead, for the median probability model. Frequently, however, the median probability model will coincide with the highest posterior probability model. One obvious situation in which this is the case is when there is a model with posterior probability greater than 1/2. Indeed, when the highest posterior probability model has substantially larger probability than the other models, it will typically also be the median probability model. Another situation in which the two coincide is when

$$(10) \qquad P(M_{\mathbf{l}}|\mathbf{y}) = \prod_{i=1}^{k} p_i^{l_i}(1 - p_i)^{(1-l_i)},$$

where the $p_i$ are the posterior inclusion probabilities in (6). This will be seen to occur in the problem of variable selection under an orthogonal design matrix, certain prior structures, and known variance $\sigma^2$, as in George and McCulloch (1993). [Clyde and Parmigiani (1996) and Clyde, DeSimone and Parmigiani (1996) show that (10) can often be approximately satisfied when $\sigma^2$ is unknown, and it is likely that the median probability model will equal the maximum probability model in such cases.]

That the median probability model is optimal (under restrictions) for both the variable selection problem and the nested case, which are very different in nature, suggests that it might quite generally be the optimal predictive model and should replace the highest posterior probability model as the "preferred" predictive model in practice. (We will see evidence of this later.) Note also that determination of the median probability model is very straightforward within ordinary MCMC model search schemes. In these schemes one develops a Markov chain to move between the models, with the posterior probability of a model being estimated by the fraction of the time that the model is visited by the chain. To determine the median probability model one need only record, for each variable, the fraction of the time that the variable is present in the visited models; at the end of the MCMC one chooses those variables for which the fraction exceeds $1/2$. Indeed, determining the median probability model in this fashion will often be computationally simpler than finding the highest posterior probability model. In the variable selection problem, for instance, accurately determining when $k$ fractions are above or below $1/2$ is often much easier than trying to accurately estimate the fractional visiting times of $2^k$ models. Note also that in the orthogonal design situation mentioned above the posterior inclusion probabilities are actually available in closed form.

The difference between predictive optimality and highest posterior model probability also explains several misunderstandings that have arisen out of the literature. For instance, Shibata (1983) shows that the BIC model selection criterion is asymptotically inferior to AIC for prediction in scenarios such as polynomial regression, when the true regression function is not a polynomial. This has been routinely misinterpreted as saying that the Bayesian approach to model selection is fine if the true model is among those being considered, but is inferior if the true model is outside the set of candidate models. Note, however, that BIC is essentially just an approximation to the log posterior probability of a model, so that model selection according to BIC is (at best) just selecting the highest posterior probability model, which is often not the optimal Bayesian answer. Indeed, as discussed above, the optimal Bayesian predictive model in the situation of Shibata (1983) is actually the median probability model. [There are also concerns with the applicability of BIC as an approximation to log posterior probability here; see Berger, Ghosh and Mukhopadhyay (2003) for further discussion.]

In Section 2, we set the basic notation for the prediction problem and give the formula for predictive expected loss. Section 3 gives the basic theory concerning optimality of the median probability model, and discusses application in nested model and ANOVA scenarios. Section 4 generalizes the basic theory to deal with problems in which all models have common nuisance parameters and the design matrix is nonorthogonal. A geometric description of the problem is provided in Section 5; this provides considerable insight into the structure of the problem. Finally, Section 6 gives some concluding comments, primarily relating to the limitations of the theory.

## 2. Preliminaries.

2.1. *Posterior inputs to the prediction problem.* Information from the data and prior is summarized by providing, for all $\mathbf{l}$,

$$
\begin{aligned}
(11) \qquad & p_{\mathbf{l}} \equiv P(M_{\mathbf{l}}|\mathbf{y}) && \text{the posterior probability of model } M_{\mathbf{l}}, \\
& \pi_{\mathbf{l}}(\boldsymbol{\beta}_{\mathbf{l}}, \sigma|\mathbf{y}) && \text{the posterior distribution} \\
& && \text{of the unknown parameters in } M_{\mathbf{l}}.
\end{aligned}
$$

These inputs will often arise from a pure Bayesian analysis based on initial specification of prior probabilities $P(M_{\mathbf{l}})$ for the models, together with prior distributions $\pi_{\mathbf{l}}(\boldsymbol{\beta}_{\mathbf{l}}, \sigma)$ for the corresponding unknown parameters. Then, given the data $\mathbf{y}$, the posterior probability of $M_{\mathbf{l}}$ is given by

$$
(12) \qquad p_{\mathbf{l}} = \frac{P(M_{\mathbf{l}})m_{\mathbf{l}}(\mathbf{y})}{\sum_{\mathbf{l}^*} P(M_{\mathbf{l}^*})m_{\mathbf{l}^*}(\mathbf{y})},
$$

where

$$
(13) \qquad m_{\mathbf{l}}(\mathbf{y}) = \int \pi_{\mathbf{l}}(\boldsymbol{\beta}_{\mathbf{l}}, \sigma) f_{\mathbf{l}}(\mathbf{y}|\boldsymbol{\beta}_{\mathbf{l}}, \sigma) \, d\boldsymbol{\beta}_{\mathbf{l}} \, d\sigma
$$

is the marginal density of $\mathbf{y}$ under $M_{\mathbf{l}}$, with $f_{\mathbf{l}}(\mathbf{y}|\boldsymbol{\beta}_{\mathbf{l}}, \sigma)$ denoting the normal density specified by $M_{\mathbf{l}}$. Likewise, the posterior distributions $\pi_{\mathbf{l}}(\boldsymbol{\beta}_{\mathbf{l}}, \sigma|\mathbf{y})$ are given by straightforward application of Bayes theorem within each model.

We allow, however, for nontraditional determination of the $p_{\mathbf{l}}$ and $\pi_{\mathbf{l}}(\boldsymbol{\beta}_{\mathbf{l}}, \sigma|\mathbf{y})$, as can occur with use of default strategies. In particular, it is not uncommon to use separate methodologies to arrive at the $p_{\mathbf{l}}$ and the $\pi_{\mathbf{l}}(\boldsymbol{\beta}_{\mathbf{l}}, \sigma|\mathbf{y})$, the $p_{\mathbf{l}}$ being determined through use of a default model selection tool such as BIC, Intrinsic Bayes Factors [cf. Berger and Pericchi (1996a)] or Fractional Bayes Factors [cf. O'Hagan (1995)]; and the $\pi_{\mathbf{l}}(\boldsymbol{\beta}_{\mathbf{l}}, \sigma|\mathbf{y})$ being determined from ordinary noninformative priors, typically the reference priors, which are either constant in the known variance case or given by

$$
(14) \qquad \pi_{\mathbf{l}}(\boldsymbol{\beta}_{\mathbf{l}}, \sigma) = \frac{1}{\sigma}
$$

in the unknown variance case. Of course, this may be incoherent from a Bayesian perspective, since it essentially means using different priors to determine the $p_{\mathbf{l}}$ and the $\pi_{\mathbf{l}}(\boldsymbol{\beta}_{\mathbf{l}}, \sigma|\mathbf{y})$.

The result of following such a mixed strategy also allows non-Bayesians to connect with this methodology. In particular, the predictor that results from use of the reference prior in model $M_\mathbf{1}$ is easily seen to be the usual least squares predictor, based on the least squares estimate

$$\widehat{\boldsymbol{\beta}}_\mathbf{1} = (\mathbf{X}_\mathbf{1}'\mathbf{X}_\mathbf{1})^{-1}\mathbf{X}_\mathbf{1}'\mathbf{y}.$$

Thus, for instance, the common use of BIC together with least squares estimates can be converted into our setting by defining $p_\mathbf{1} \propto e^{\mathrm{BIC}}$.

Finally, the empirical Bayes approach is often used to obtain estimated versions of the quantities in (11). Again, while not strictly coherent from a Bayesian perspective, one can utilize such inputs in the following methodology.

2.2. *Predictors and predictive expected loss.* It is easy to see that the optimal predictor of $y^*$ in (3), under squared error loss and when the model $M_\mathbf{1}$, of dimension $k_l$, is true, is given by

$$\widehat{y}_\mathbf{1}^* = \mathbf{x}^*\mathbf{H}_\mathbf{1}\widetilde{\boldsymbol{\beta}}_\mathbf{1},$$

where $\widetilde{\boldsymbol{\beta}}_\mathbf{1}$ is the posterior mean of $\boldsymbol{\beta}_\mathbf{1}$ with respect to $\pi_\mathbf{1}(\boldsymbol{\beta}_\mathbf{1}, \sigma | \mathbf{y})$ and $\mathbf{H}_\mathbf{1}$ is the $k \times k_l$ matrix whose $(i,j)$ entry is 1 if $l_i = 1$ and $j = \sum_{r=1}^{i} l_r$ and is 0 otherwise. Note that $\mathbf{H}_\mathbf{1}$ is simply the matrix such that $\mathbf{x}\mathbf{H}_\mathbf{1}$ is the subvector of $\mathbf{x}$ corresponding to the nonzero coordinates of $\mathbf{l}$, that is, the covariate vector corresponding to model $M_\mathbf{1}$. The posterior mean of $\boldsymbol{\beta}$ in the full model is thus formally written as $\widetilde{\boldsymbol{\beta}}_{(1,\ldots,1)}$, but we will drop the subscript and simply denote it by $\widetilde{\boldsymbol{\beta}}$ (as we have done with $\widehat{\boldsymbol{\beta}}$, the least squares estimate for the full model).

The optimal Bayesian predictor of $y^*$ is well known to be the model averaging predictor, given by

$$(15) \qquad \bar{y}^* = \mathbf{x}^*\bar{\boldsymbol{\beta}} \equiv \mathbf{x}^* \sum_\mathbf{1} p_\mathbf{1}\mathbf{H}_\mathbf{1}\widetilde{\boldsymbol{\beta}}_\mathbf{1}.$$

Note that "optimal" is defined in terms of expected loss over the posterior in (11); if the posterior arises from an incoherent prior, as is allowed above, there is no guarantee that the resulting procedure possesses any other optimality properties—even "dutch books" against the procedure or frequentist inadmissibility could result.

The best single model for prediction can be found by the following lemma.

LEMMA 1. *The optimal model for prediction of $y^*$ in (3) under the squared error loss (4), when the future covariates satisfy (5) and the posterior distribution is as in (11), is the model that minimizes*

$$(16) \qquad R(M_\mathbf{1}) \equiv (\mathbf{H}_\mathbf{1}\widetilde{\boldsymbol{\beta}}_\mathbf{1} - \bar{\boldsymbol{\beta}})'\mathbf{Q}(\mathbf{H}_\mathbf{1}\widetilde{\boldsymbol{\beta}}_\mathbf{1} - \bar{\boldsymbol{\beta}}),$$

*where $\bar{\boldsymbol{\beta}}$ is defined in (15).*

PROOF. For fixed $\mathbf{x}^*$ a standard result [see, e.g., Bernardo and Smith (1994), page 398] is that

$$\mathbf{E}[(\widehat{y}_l^* - y^*)^2] = C + (\widehat{y}_l^* - \bar{y}^*)^2,$$

where $C$ does not depend on $l$ and the expectation is with respect to the predictive distribution of $y^*$ given $y$. Since

$$(\widehat{y}_l^* - \bar{y}^*)^2 = (\mathbf{H}_l \widetilde{\boldsymbol{\beta}}_l - \bar{\boldsymbol{\beta}})' \mathbf{x}^{*\prime} \mathbf{x}^* (\mathbf{H}_l \widetilde{\boldsymbol{\beta}}_l - \bar{\boldsymbol{\beta}}),$$

taking the expectation over $\mathbf{x}^*$ and using (5) yields the result.  □

**3. Basic results and examples.** Section 3.1 presents the basic theorem that is used to establish optimality of the median probability model. Section 3.2 considers the situation in which all submodels of the linear model are allowed. Section 3.3 deals with nested models and Section 3.4 considers the ANOVA situation.

3.1. *Basic theory.* Assume $\mathbf{X}'\mathbf{X}$ is diagonal. Then it will frequently be the case that the posterior means $\widetilde{\boldsymbol{\beta}}_l$ satisfy

(17)
$$\widetilde{\boldsymbol{\beta}}_l = \mathbf{H}_l'\widetilde{\boldsymbol{\beta}},$$

that is, that the posterior mean of $\boldsymbol{\beta}_l$ is found by simply taking the relevant coordinates of $\widetilde{\boldsymbol{\beta}}$, the posterior mean in the full model. Here are two common scenarios in which (17) is true.

CASE 1. *Noninformative priors for model parameters*: Use of the reference priors in (14) (or constant priors when $\sigma^2$ is known) results in the posterior means being the least squares estimates $\widehat{\boldsymbol{\beta}}_l$. Because $\mathbf{X}'\mathbf{X}$ is diagonal, it is easy to see that (17) is then satisfied.

CASE 2. *Independent conjugate normal priors*: In the full model suppose that $\pi(\boldsymbol{\beta}|\sigma)$ is $\mathcal{N}_k(\boldsymbol{\mu}, \sigma^2\boldsymbol{\Lambda})$, the $k$-variate normal distribution with mean $\boldsymbol{\mu}$ and diagonal covariance matrix $\sigma^2\boldsymbol{\Lambda}$, with $\boldsymbol{\Lambda}$ given. Then it is natural to choose the priors on $\boldsymbol{\beta}_l$ in the submodels to be $\mathcal{N}_{k_l}(\mathbf{H}_l'\boldsymbol{\mu}, \sigma^2\mathbf{H}_l'\boldsymbol{\Lambda}\mathbf{H}_l)$, where $k_l$ is the dimension of $\boldsymbol{\beta}_l$. It is then easy to verify that (17) holds for any prior on $\sigma^2$ or for $\sigma^2$ being given (e.g., known, or estimated). Note that we do not necessarily recommend using this conjugate form of the prior with unknown $\sigma^2$; see Berger and Pericchi (2001) for discussion.

While $\boldsymbol{\Lambda}$ could be chosen subjectively, it is more common to utilize default choices, such as the *g-type* normal priors [cf. Zellner (1986)] $\boldsymbol{\Lambda} = n(\mathbf{X}'\mathbf{X})^{-1}$ or $\boldsymbol{\Lambda} = c(\mathbf{X}'\mathbf{X})^{-1}$, with $c$ chosen by an empirical Bayes analysis (e.g., chosen to maximize the marginal density averaged over models). Papers which fall under these settings include Chipman, George and McCulloch (2001), Clyde and Parmigiani (1996), Clyde, DeSimone and Parmigiani (1996), Clyde, Parmigiani and Vidakovich (1998) and George and Foster (2000).

Note that one can use noninformative priors for certain coordinates and independent conjugate normal priors for other coordinates. This is particularly useful when all models under consideration have "common" unknown parameters; it is then typical to utilize noninformative priors for the common parameters, while using independent conjugate normal priors for the other parameters [see Berger, Pericchi and Varshavsky (1998) for justification of this practice].

LEMMA 2. *If* $\mathbf{Q}$ *is diagonal with diagonal elements* $q_i > 0$ *and* (17) *holds, then*

$$(18) \qquad R(M_{\mathbf{l}}) = \sum_{i=1}^{k} \widetilde{\beta}_i^2 q_i (l_i - p_i)^2,$$

*where* $p_i$ *is as in* (6).

PROOF. From (17) it follows that

$$\bar{\boldsymbol{\beta}} \equiv \sum_{\mathbf{l}} p_{\mathbf{l}} \mathbf{H}_{\mathbf{l}} \widetilde{\boldsymbol{\beta}}_{\mathbf{l}} = \sum_{\mathbf{l}} p_{\mathbf{l}} \mathbf{H}_{\mathbf{l}} \mathbf{H}_{\mathbf{l}}' \widetilde{\boldsymbol{\beta}} = \mathbf{D}(\mathbf{p}) \widetilde{\boldsymbol{\beta}},$$

where $\mathbf{D}(\mathbf{p})$ is the diagonal matrix with diagonal elements $p_i$. Likewise, (16) becomes

$$R(M_{\mathbf{l}}) = (\mathbf{H}_{\mathbf{l}} \mathbf{H}_{\mathbf{l}}' \widetilde{\boldsymbol{\beta}} - \mathbf{D}(\mathbf{p}) \widetilde{\boldsymbol{\beta}})' \mathbf{Q} (\mathbf{H}_{\mathbf{l}} \mathbf{H}_{\mathbf{l}}' \widetilde{\boldsymbol{\beta}} - \mathbf{D}(\mathbf{p}) \widetilde{\boldsymbol{\beta}})$$

$$= \widetilde{\boldsymbol{\beta}}' (\mathbf{D}(\mathbf{l}) - \mathbf{D}(\mathbf{p})) \mathbf{Q} (\mathbf{D}(\mathbf{l}) - \mathbf{D}(\mathbf{p})) \widetilde{\boldsymbol{\beta}},$$

and the conclusion is immediate. □

THEOREM 1. *If* $\mathbf{Q}$ *is diagonal with diagonal elements* $q_i > 0$, *condition* (17) *holds and the models under consideration have graphical model structure, then the median probability model is the best predictive model.*

PROOF. To minimize (18) among all possible models, it is clear that one should choose $l_i = 1$ if $p_i \geq 1/2$ and $l_i = 0$ otherwise, which is as in (7). As mentioned earlier, the graphical model structure ensures that the model so defined is actually in the space of models under consideration, completing the proof. □

The above theorem did not formally use the condition that $\mathbf{X}'\mathbf{X}$ be diagonal. However, if it is not diagonal, then (17) will not typically hold, nor will $\mathbf{Q}$ usually be diagonal.

3.2. *All submodels scenario.* Under the same conditions as in Section 3.1 the following corollary to Theorem 1 gives the median probability model when all submodels are considered.

COROLLARY 1. *If $\mathbf{Q}$ is diagonal with diagonal elements $q_i > 0$, condition (17) holds and any submodel of the full model is allowed, then the best predictive model is the median probability model given by (7). In addition, if $\sigma^2$ is given in Case 2 of Section 3.1 and the prior probabilities of the models satisfy*

$$P(M_\mathbf{l}) = \prod_{i=1}^{k} (p_i^0)^{l_i} (1 - p_i^0)^{(1-l_i)}, \tag{19}$$

*where $p_i^0$ is the prior probability that variable $x_i$ is in the model, then (10) is satisfied and the median probability model is the model with highest posterior probability.*

PROOF. The first part of the corollary is immediate, since {all submodels} clearly has graphical model structure.

For Case 2 and given $\sigma^2$, computation as in Clyde and Parmigiani (1996) and Clyde, DeSimone and Parmigiani (1996) shows that (10) is satisfied with

$$p_i^{-1} = 1 + \left( \frac{1}{p_i^0} - 1 \right) (1 + \lambda_i d_i)^{1/2} \exp \left\{ -\frac{v_i^2 \lambda_i + 2v_i \mu_i - \mu_i^2 d_i}{2\sigma^2 (1 + \lambda_i d_i)} \right\},$$

where the $\{d_i\}$ and $\{\lambda_i\}$ are the diagonal elements of $\mathbf{X}'\mathbf{X}$ and $\mathbf{\Lambda}$, respectively, and $\mathbf{v} = (v_1, \ldots, v_k)' = \mathbf{X}'\mathbf{y}$. □

While many proposed choices of prior probabilities satisfy (19), others do not. For instance, Jeffreys (1961) suggested that it might often be reasonable to choose the prior probability of given model orders to be decreasing, for example, $P(\text{order } j) \propto 1/j$, with this probability then being divided up equally among all models of size $j$. Such an assignment of prior probabilities would not typically satisfy (19), and the best predictive model (i.e., the median probability model) would then not necessarily be the highest posterior probability model, even in Case 2.

Finally, it should be noted that Corollary 1 also applies to the case where all models under consideration have "common parameters." One can simply define $p_i^0 = 1$ for such parameters.

3.3. *Nested models.* We initially consider the orthogonal case, as in Section 3.1. This is generalized to the nonorthogonal case at the end of the section.

### 3.3.1. *Orthogonal case.*

COROLLARY 2. *If $\mathbf{Q}$ is diagonal with diagonal elements $q_i > 0$, condition ($17$) holds and the models under consideration are nested, then the best predictive model is the median probability model given by ($7$) or ($9$).*

PROOF. This is immediate from Theorem [3.1], because nested models have graphical model structure. $\square$

EXAMPLE 1. *Nonparametric regression* [also studied in Mukhopadhyay ($2000$) for a different purpose]: The data consists of the paired observations $(x_i, y_i), i = 1, \ldots, n$, where for known $\sigma^2$,

$$(20) \qquad y_i = f(x_i) + \varepsilon_i, \qquad \varepsilon_i \sim N(0, \sigma^2).$$

Represent $f(\cdot)$, an unknown function defined on the interval $(-1, 1)$, using an orthonormal series expansion, as

$$f(x) = \sum_{i=1}^{\infty} \beta_i \phi_i(x),$$

where $\{\phi_1(x), \phi_2(x), \ldots\}$ are the Chebyshev polynomials. Of course, only a finite approximation to this series can be utilized, so define the model $M_{\mathbf{l}(j)}$ [see ($8$) for notation] to be

$$M_{\mathbf{l}(j)} : y = \sum_{i=1}^{j} \beta_i \phi_i(x) + \varepsilon, \qquad \varepsilon \sim N(0, \sigma^2).$$

The problem is thus to choose among the nested sequence of linear models $M_{\mathbf{l}(j)}, j = 1, \ldots$. As is common in practice, we will choose an upper bound $k$ on the size of the model, so that $M_{\mathbf{l}(k)}$ is the full model in our earlier terminology.

The function $f(x) = -\log(1 - x)$ was considered in Shibata ($1983$) as an example for which BIC yielded considerably suboptimal models for prediction. It is hence of interest to see how the median probability model fares in this situation.

We assume the $y_i$ are observed at the covariates $x_i = \text{Cosine}([n - i + \frac{1}{2}]\frac{\pi}{n})$, $i = 1, \ldots, n$, and let $\mathbf{X}_j = (\phi_m(x_i))$ be the resulting $n \times j$ design matrix with indicated $(i, m)$ entries, $i = 1, \ldots, n$ and $m = 1, \ldots, j$. From the definition of Chebyshev polynomials it follows that $\mathbf{X}_j' \mathbf{X}_j = \frac{n}{2} \mathbf{I}_j$, where $\mathbf{I}_j$ is the $j \times j$ identity matrix. It follows that the least squares estimate of $\boldsymbol{\beta}_j$ is $\widehat{\boldsymbol{\beta}}_j = \frac{2}{n} \mathbf{X}_j' \mathbf{y}$.

Assume that the models have equal prior probability $1/k$, and that, within any model $M_{\mathbf{l}(j)}$ the $\beta_i$ have independent $N(0, ci^{-a})$ prior distributions for some constants $c$ and $a$ (which are the same across models). The choice of $a$ determines how quickly the prior variances of the $\beta_i$ decrease (any $L_2$ function must have $a > 1$), and we shall consider three choices: $a = 1$, $a = 2$ (which happens to be the rate corresponding to the test function) and $a = 3$. For simplicity of calculation we estimate $c$ by an empirical Bayes analysis using the full model $M_{I(k)}$, keeping the estimate $\hat{c}$ fixed across models. Then if $\mathbf{Q}$ is diagonal (as would thus be the case for the natural choice $\mathbf{Q} = \mathbf{X}_k' \mathbf{X}_k$), Corollary 2 implies that the median probability model will be the optimal Bayesian predictive model.

For nonparametric regression it is common to utilize the loss function

$$(21) \qquad L(f, \hat{f}) = \int_{-1}^{1} (\hat{f}(x) - f(x))^2 \, dx.$$

In the predictive context use of this loss is equivalent to prediction under squared error loss when the future covariates $x$ are thought to be uniformly distributed on $(-1, 1)$. A standard computation shows that $L(f, \hat{f}) = \sum_{i=1}^{\infty} (\widetilde{\beta}_i - \beta_i)^2$, where $\widetilde{\beta}_i$ stands for the estimator that is used for the true coefficient $\beta_i$. Since we have restricted the models under consideration to be of maximum order $k$, it follows that $\widetilde{\beta}_i = 0$ for $i > k$ in our problem and the loss arising from these coordinates can be ignored in model comparison. The resulting loss is also easily seen to be equivalent to the predictive loss we have previously used, with $\mathbf{Q} = \mathbf{I}_k$.

The optimality of the median probability model is with respect to the internal Bayesian computation, which assumes that the true function is actually one of the models $M_{I(j)}$, $j = 1, \ldots, k$. For the example considered by Shibata, however, the true model lies outside the collection of Bayesian models (since it can be shown that none of the $\beta_i$ in the expansion for this function is zero). It is thus of interest to see how the median probability model performs in terms of the loss (21) for the (here known) function.

Under $M_{I(j)}$ the estimates of the $\beta_i$ are the (empirical) Bayes estimates $\widetilde{\beta}_i = (1 + 2\sigma^2 i^a / n\hat{c})^{-1} \hat{\beta}_i$ if $i \leq k$, and $\widetilde{\beta}_i = 0$ otherwise. Hence the predictive loss under model $M_{I(j)}$ is $\sum_{i=1}^{j} (\widetilde{\beta}_i - \beta_i)^2 + \sum_{i=j+1}^{\infty} \beta_i^2$, although we will ignore the terms in the sum for $i > k$ since they are common across all considered models.

In Table 1 we compare the expected predictive loss (the frequentist expectation with respect to the data under the true function) of the maximum probability model with that of the median probability model. We also include the model averaged estimate arising from (15) in the comparison; this is the optimal estimate from the internal Bayesian perspective. Finally, we also consider AIC and BIC. These model selection criteria are most commonly used in conjunction with least squares parameter estimates, and that choice is made for computing the expected predictive losses in Table 1.

The entries in the table were computed by simulating data from the true function, selecting a model, computing the corresponding function estimate, and finally determining the actual loss. This was repeated a total of $N = 1000, 1000, 100$ times, respectively, for the three cases in the table, with the resulting averages forming the table entries. Note that the largest model sizes considered for the three cases were $k = 29, 79, 79$, respectively.

Our main goal was to compare the maximum probability model and the median probability model. The median probability model is clearly significantly better, even in terms of this frequentist expected loss. (Again, we know it is better in terms of Bayesian predictive loss.) Indeed, the median probability model is almost as good as the model averaging estimate (and in two cases is even better); since model averaging is thought of as optimal, its near equivalence with MedianProb is compelling.

Note that AIC does seem to do better than BIC, as reported by Shibata, but all of the actual Bayesian procedures are considerably better than either. This is in part due to the fact that the Bayesian procedures use better parameter estimates than least squares, but the dominance of MedianProb and ModelAv (but not MaxProb) can be seen to hold even if the superior shrinkage estimates are also used with BIC and AIC.

3.3.2. *Nonorthogonal case.* Corollary 2 presented a result for nested models in the case of an orthogonal design matrix. It is not too surprising that

TABLE 1
*For various $n$ and $\sigma^2$ the expected loss and average model size for the maximum probability model (MaxProb), the median probability model (MedianProb), model averaging (ModelAv), and BIC and AIC, in the Shibata example*

| | Expected loss [average model size] | | | | |
|---|---|---|---|---|---|
| | **MaxProb** | **MedianProb** | **ModelAv** | **BIC** | **AIC** |
| $n = 30, \sigma^2 = 1$ | | | | | |
| $a = 1$ | 0.99 [8] | 0.89 [10] | 0.84 | 1.14 [8] | 1.09 [7] |
| $a = 2$ | 0.88 [10] | 0.80 [16] | 0.81 | 1.14 [8] | 1.09 [7] |
| $a = 3$ | 0.88 [9] | 0.84 [17] | 0.85 | 1.14 [8] | 1.09 [7] |
| $n = 100, \sigma^2 = 1$ | | | | | |
| $a = 1$ | 0.54 [14] | 0.51 [19] | 0.47 | 0.59 [11] | 0.59 [13] |
| $a = 2$ | 0.47 [23] | 0.43 [43] | 0.44 | 0.59 [11] | 0.59 [13] |
| $a = 3$ | 0.47 [22] | 0.46 [45] | 0.46 | 0.59 [11] | 0.59 [13] |
| $n = 2000, \sigma^2 = 3$ | | | | | |
| $a = 1$ | 0.34 [23] | 0.33 [26] | 0.30 | 0.41 [12] | 0.38 [21] |
| $a = 2$ | 0.26 [42] | 0.25 [51] | 0.25 | 0.41 [12] | 0.38 [21] |
| $a = 3$ | 0.29 [38] | 0.29 [50] | 0.29 | 0.41 [12] | 0.38 [21] |

such results can be obtained under orthogonality. Quite surprising, however, is that the orthogonality condition can be removed under the following two conditions:

CONDITION 1.    $\mathbf{Q} = \gamma \mathbf{X}' \mathbf{X}$ for some $\gamma > 0$.

CONDITION 2.    $\widetilde{\boldsymbol{\beta}}_\mathbf{l} = b \widehat{\boldsymbol{\beta}}_\mathbf{l}$, where $b > 0$, that is, the posterior means are proportional to the least squares estimates, with the same proportionality constant across models.

Note that Condition 2 is merely a special case of (17), so the generality here is in Condition 1, allowing for nondiagonal $\mathbf{Q}$. (Of course, $\mathbf{Q} = \gamma \mathbf{X}' \mathbf{X}$ is diagonal under orthogonality.)

There are two common scenarios in which Condition 2 is satisfied. The first is when the reference priors (14) are used, in which case the posterior means are the least squares estimates. The second is when using *g-type* normal priors [cf. Zellner (1986)], where $\pi_\mathbf{l}(\boldsymbol{\beta}_\mathbf{l} | \sigma)$ is $\mathcal{N}_{k_\mathbf{l}}(\mathbf{0}, c\sigma^2 (\mathbf{X}_\mathbf{l}' \mathbf{X}_\mathbf{l})^{-1})$, with the same constant $c > 0$ for each model. (This constant could be specified or estimated in an empirical Bayesian context.) It is then easy to verify that Condition 2 holds with $b = c/(1 + c)$ (irrespective of the prior for $\sigma$).

THEOREM 2.    *For a sequence of nested models for which Conditions* 1 *and* 2 *hold, the best predictive model is the median probability model given by* (7) *or* (9).

PROOF.    From (16) and using Conditions 1 and 2,

$$(22) \qquad R(M_{\mathbf{l}(j)}) = \gamma b^2 (\mathbf{H}_{\mathbf{l}(j)} \widehat{\boldsymbol{\beta}}_{\mathbf{l}(j)} - \bar{\boldsymbol{\beta}})' \mathbf{X}' \mathbf{X} (\mathbf{H}_{\mathbf{l}(j)} \widehat{\boldsymbol{\beta}}_{\mathbf{l}(j)} - \bar{\boldsymbol{\beta}}).$$

Noting that $\mathbf{X} \mathbf{H}_{\mathbf{l}(j)} = \mathbf{X}_{\mathbf{l}(j)}$ and defining $\mathbf{P}_\mathbf{l} = \mathbf{X}_\mathbf{l} (\mathbf{X}_\mathbf{l}' \mathbf{X}_\mathbf{l})^{-1} \mathbf{X}_\mathbf{l}'$, it follows that

$$R(M_{\mathbf{l}(j)}) = \gamma b^2 \mathbf{y}' \left( \mathbf{P}_{\mathbf{l}(j)} - \sum_{i=1}^{k} p_{\mathbf{l}(i)} \mathbf{P}_{\mathbf{l}(i)} \right)^2 \mathbf{y}.$$

Note that $\mathbf{P}_\mathbf{l}^2 = \mathbf{P}_\mathbf{l}$ and $\mathbf{P}_{\mathbf{l}(i)} \mathbf{P}_{\mathbf{l}(j)} = \mathbf{P}_{\mathbf{l}(\min\{i,j\})}$. Hence, expanding the quadratic in (22) yields

$$R(M_{\mathbf{l}(j)}) = \gamma b^2 \mathbf{y}' \left( \sum_{i=1}^{k} p_{\mathbf{l}(i)} \mathbf{P}_{\mathbf{l}(i)} \right)^2 \mathbf{y}$$
$$+ \gamma b^2 \mathbf{y}' \left( \mathbf{P}_{\mathbf{l}(j)} - 2 \sum_{i=1}^{j-1} p_{\mathbf{l}(i)} \mathbf{P}_{\mathbf{l}(i)} - 2 \sum_{i=j}^{k} p_{\mathbf{l}(i)} \mathbf{P}_{\mathbf{l}(j)} \right) \mathbf{y}.$$

It follows that

$$R(M_{\mathbf{l}(j+1)}) - R(M_{\mathbf{l}(j)}) = \gamma b^2 \left(1 - 2\sum_{i=j+1}^{k} p_{\mathbf{l}(i)}\right) \mathbf{y}'(\mathbf{P}_{\mathbf{l}(j+1)} - \mathbf{P}_{\mathbf{l}(j)})\mathbf{y}.$$

Since $\mathbf{y}'(\mathbf{P}_{\mathbf{l}(j+1)} - \mathbf{P}_{\mathbf{l}(j)})\mathbf{y} > 0$ and the $(1 - 2\sum_{i=j+1}^{k} p_{\mathbf{l}(i)})$ are increasing in $j$ from $-1$ to $+1$, moving to a larger model will reduce the risk until $(1 - 2\sum_{i=j+1}^{k} p_{\mathbf{l}(i)})$ first turns positive. The conclusion is immediate. $\square$

An example of a nested model in the nonorthogonal case will be given in Section 4.

3.4. *ANOVA.* Many ANOVA problems, when written in linear model form, yield diagonal $\mathbf{X}'\mathbf{X}$ and any such problems will naturally fit under the theory of Section 3.1. In particular, this is true for any balanced ANOVA in which each factor has only two levels.

To see the idea, it suffices to consider the case of two factors A and B each with two levels. The full two-way ANOVA model with interactions is

$$y_{ijk} = \mu + a_i + b_j + ab_{ij} + \varepsilon_{ijk}$$

with $i = 1, 2$, $j = 1, 2$, $k = 1, 2, \ldots, K$ and $\varepsilon_{ijk}$ independent $N(0, \sigma^2)$, with $\sigma^2$ unknown. In our earlier notation, this can be written

$$\mathbf{y} = \mathbf{X}\boldsymbol{\beta} + \boldsymbol{\varepsilon},$$

where

$$\mathbf{y} = (y_{111}, \ldots, y_{11K}, y_{121}, \ldots, y_{12K}, y_{211}, \ldots, y_{21K}, y_{221}, \ldots, y_{22K})',$$

$$\boldsymbol{\beta} = (\mu, a_1, b_1, ab_{11})'$$

and $\mathbf{X}$ is the $4K \times 4$ matrix

$$X = \begin{pmatrix} 1 & 1 & 1 & 1 \\ \vdots & \vdots & \vdots & \vdots \\ 1 & 1 & 1 & 1 \\ 1 & 1 & -1 & -1 \\ \vdots & \vdots & \vdots & \vdots \\ 1 & 1 & -1 & -1 \\ 1 & -1 & 1 & -1 \\ \vdots & \vdots & \vdots & \vdots \\ 1 & -1 & 1 & -1 \\ 1 & -1 & -1 & 1 \\ \vdots & \vdots & \vdots & \vdots \\ 1 & -1 & -1 & 1 \end{pmatrix},$$

where the last column is the product of the second and the third, since $a_1 = -a_2$, $b_1 = -b_2$, $ab_{11} = ab_{22} = -ab_{12} = -ab_{21}$. Computation then shows that $\mathbf{X}'\mathbf{X} = 4K\mathbf{I}_4$, so that the earlier theory will apply.

There are several model comparison scenarios of interest. We use a slight modification of the previous model notation for simplicity, for example, $M_{1011}$ instead of $M_{(1,0,1,1)}$, representing the model having all parameters except $a_1$.

SCENARIO 1. *All models with the constant* $\mu$: Thus the set of models under consideration is

$$\{M_{1000}, M_{1100}, M_{1010}, M_{1001}, M_{1101}, M_{1011}, M_{1110}, M_{1111}\}.$$

SCENARIO 2. *Interactions present only with main effects, and* $\mu$ *included*: The set of models under consideration here is $\{M_{1000}, M_{1100}, M_{1010}, M_{1110}, M_{1111}\}$. Note that this set of models has graphical structure.

SCENARIO 3. *An analogue of an unusual classical test*: In classical ANOVA testing it is sometimes argued [cf. Scheffé (1959), pages 94 and 110] that one might be interested in testing for no interaction effect followed by testing for the main effects, even if the no-interaction test rejected. (It is argued that the hypotheses of zero main effects could still be accepted, which would imply that, while there are differences, the tests do not demonstrate any differences in the levels of one factor averaged over the levels of the other.) The four models that are under consideration in this process, including the constant $\mu$ in all, are $\{M_{1101}, M_{1011}, M_{1110}, M_{1111}\}$.

We do not comment upon the reasonableness of considering this class of models, but are interested in the class because it does *not* have graphical model structure and yet the median probability model is guaranteed to be in the class. To see this, consider the possibility that $a_1$ has posterior inclusion probability less than $1/2$, and so would be excluded from the median probability model. Clearly this can only happen if $M_{1011}$ has posterior probability greater than $1/2$; but then $M_{1011}$ would automatically be the median probability model. Arguing similarly for $b_1$ and $ab_{11}$, one can conclude that $M_{1111}$ will be the median probability model, unless one of the other models has posterior probability greater than $1/2$, in which case it will be the median probability model.

EXAMPLE 2. Montgomery [(1991), pages 271–274] considers the effects of the concentration of a reactant and the amount of a catalyst on the yield in a chemical process. The reactant concentration is factor A and has two levels, 15% and 25%. The catalyst is factor B, with the two levels "one bag" and "two bags" of catalyst. The experiment was replicated three times and the data are given in Table 2. Note that the classical ANOVA tests of "no A effect," "no B effect" and "no interaction effect" resulted in $p$-values of 0.00008, 0.00236 and 0.182, respectively. The Bayesian quantities that would be used analogously are the posterior variable inclusion probabilities, $p_2$, $p_3$ and $p_4$. (Strictly, $1 - p_i$ would be the analogue of the corresponding $p$-value.)

To carry out the Bayesian analysis, the reference prior $\pi(\mu, \sigma) \propto \frac{1}{\sigma}$ was used for the common parameters, while the standard $N(0, \sigma^2)$ $g$-prior was used for $a_1$, $b_1$ and $ab_{11}$. In each scenario the models under consideration were given equal prior probabilities of being true. The conditions of Section 3.1 are then satisfied, so that we know that the median probability model will be the optimal predictive model. For the three scenarios described above the results of the Bayesian analysis are given in Tables 3–5. In all three scenarios the median probability model indeed has the lowest posterior expected loss (as was known from the theory). Interestingly, the median probability model equals the maximum probability model in all three scenarios and is the model $M_{1110}$. The variable inclusion probabilities show clearly that an "A

TABLE 2
*Data for the $2^2$ ANOVA example*

| Treatment combination | Replicates | | |
|---|---|---|---|
| | **I** | **II** | **III** |
| A low, B low | 28 | 25 | 27 |
| A high, B low | 36 | 32 | 32 |
| A low, B high | 18 | 19 | 23 |
| A high, B high | 31 | 30 | 29 |

TABLE 3
*Scenario 1 (all models). Posterior probabilities and expected losses for the models. The posterior inclusion probabilities are $p_2 = 0.9977$, $p_3 = 0.9514$ and $p_4 = 0.3621$; thus $M_{1110}$ is the median probability model*

| Model | Posterior probability | Posterior expected loss |
|---|---|---|
| $M_{1000}$ | 0.0008 | 235.47 |
| $M_{1100}$ | 0.0342 | 58.78 |
| $M_{1010}$ | 0.0009 | 177.78 |
| $M_{1001}$ | 0.0003 | 237.43 |
| $M_{1110}$ | 0.6019 | 1.083 |
| $M_{1101}$ | 0.0133 | 60.73 |
| $M_{1011}$ | 0.0003 | 179.74 |
| $M_{1111}$ | 0.3483 | 3.04 |

TABLE 4

*Scenario 2 (graphical models). Posterior probabilities and expected losses for the models. The posterior inclusion probabilities are $p_2 = 0.9982$, $p_3 = 0.9644$ and $p_4 = 0.3532$; thus $M_{1110}$ is the median probability model*

| Model | Posterior probability | Posterior expected loss |
|-------|-----------------------|-------------------------|
| $M_{1000}$ | 0.0009 | 237.21 |
| $M_{1100}$ | 0.0347 | 60.33 |
| $M_{1010}$ | 0.0009 | 177.85 |
| $M_{1110}$ | 0.6103 | 0.97 |
| $M_{1111}$ | 0.3532 | 3.05 |

TABLE 5

*Scenario 3 (unusual classical models). Posterior probabilities and expected losses for the models. The posterior inclusion probabilities are $p_2 = 0.9997$, $p_3 = 0.9862$ and $p_4 = 0.3754$; thus $M_{1110}$ is the median probability model*

| Model | Posterior probability | Posterior expected loss |
|-------|-----------------------|-------------------------|
| $M_{1011}$ | 0.0003 | 180.19 |
| $M_{1101}$ | 0.0138 | 64.93 |
| $M_{1110}$ | 0.6245 | 1.01 |
| $M_{1111}$ | 0.3614 | 2.78 |

effect" and a "B effect" should be in the model (with inclusion probabilities exceeding 0.99 and 0.95, respectively), while the interaction effect has a moderately small probability (about 1/3) of being in the model.

We also carried out an analysis with the $N(0, c\sigma^2)$ $g$-prior for $a_1$, $b_1$ and $ab_{11}$, but with $c$ being estimated by maximizing the overall marginal density $\frac{1}{L} \sum_l m_l(\mathbf{y})$, where the individual marginal densities $m_l(\mathbf{y})$ are given by (13) and $L$ is the number of models under consideration. The conditions of Section 3.1 are still satisfied, so that we know that the median probability model will be the optimal predictive model. The results did not significantly change from the above analysis, however, and so are not reported.

Had $\sigma^2$ been known in Scenario 1, Corollary 1 would have established that the median probability model would equal the maximum probability model. Here $\sigma^2$ is unknown, however, and it will not always be the case that the median probability model equals the maximum probability model. Indeed, we also carried out the analysis of the example utilizing the $N(0, \sigma^2)$ $g$-prior

*Scenario 3 (unusual classical models, with g-prior for $\mu$). Posterior probabilities and expected losses for the models. The posterior inclusion probabilities are $p_2 = 0.876$, $p_3 = 0.714$ and $p_4 = 0.544$; thus $M_{1111}$ is the median probability model*

| Model | Posterior probability | Posterior expected loss |
|---|---|---|
| $M_{1011}$ | 0.124 | 143.03 |
| $M_{1101}$ | 0.286 | 36.78 |
| $M_{1110}$ | 0.456 | 10.03 |
| $M_{1111}$ | 0.134 | 9.41 |

for $\mu$, as well as for $a_1$, $b_1$ and $ab_{11}$, and found that the median probability model then differed from the maximum probability model in all three scenarios. Table 6 gives the results for Scenario 3; note that the median probability model is nearly the lowest probability model! (We do not, however, recommend this analysis; $g$-priors should not be used for parameters common to all models.)

## 4. Common nonorthogonal nuisance parameters.

Frequently all models will contain "common" parameters $\boldsymbol{\beta}_{(1)} \equiv (\beta_1, \ldots, \beta_{k_1})$. A typical example is when all models contain an overall mean $\beta_1$ [or, equivalently, when the first column of each model design matrix is $(1, \ldots, 1)'$]. For the orthogonal case discussed earlier this caused no difficulties. For the nonorthogonal case, however, this considerably complicates the analysis. Still, we will see that the median probability model remains optimal under mild modifications of the previously considered conditions.

To present the results, it is convenient to slightly change notation, writing the regression parameters of $M_\mathbf{l}$ as $(\boldsymbol{\beta}'_{(1)}, \boldsymbol{\beta}'_\mathbf{l})'$, with corresponding design matrix $(\mathbf{X}_{(1)} \mathbf{X}_\mathbf{l})$. Also, define

$$(23) \qquad \begin{aligned} \mathbf{Q}_{(1)} &= \mathbf{I} - \mathbf{X}_{(1)}(\mathbf{X}'_{(1)}\mathbf{X}_{(1)})^{-1}\mathbf{X}'_{(1)}, \\ \mathbf{P}_\mathbf{l} &= \mathbf{Q}_{(1)}^{1/2}\mathbf{X}_\mathbf{l}(\mathbf{X}'_\mathbf{l}\mathbf{Q}_{(1)}\mathbf{X}_\mathbf{l})^{-1}\mathbf{X}'_\mathbf{l}\mathbf{Q}_{(1)}^{1/2}. \end{aligned}$$

The necessary conditions are:

CONDITION 1. $\mathbf{Q} = \gamma \mathbf{X}'\mathbf{X}$ for some $\gamma > 0$.

CONDITION 2. For some fixed $b > 0$ the posterior means of $\boldsymbol{\beta}_\mathbf{l}$ are of the form $\tilde{\boldsymbol{\beta}}_\mathbf{l} = b(\mathbf{X}'_\mathbf{l}\mathbf{Q}_{(1)}\mathbf{X}_\mathbf{l})^{-1}\mathbf{X}'_\mathbf{l}\mathbf{Q}_{(1)}\mathbf{y}$.

CONDITION 3. $\pi_\mathbf{l}(\boldsymbol{\beta}_{(1)}, \boldsymbol{\beta}_\mathbf{l} | \sigma) = \pi_\mathbf{l}(\boldsymbol{\beta}_\mathbf{l} | \sigma)$ (i.e., the prior density for $\boldsymbol{\beta}_{(1)}$ in each model is constant).

A $g$-type prior for which Condition 2 holds is

$$\pi_{\mathbf{l}}(\boldsymbol{\beta}_{\mathbf{l}}|\sigma) = \mathcal{N}(\mathbf{0}, c\sigma^2(\mathbf{X}_{\mathbf{l}}'\mathbf{Q}_{(1)}\mathbf{X}_{\mathbf{l}})^{-1}),$$

with the same constant $c > 0$ for each model. It is then easy to verify that Condition 2 holds with $b = c/(1 + c)$ (irrespective of the prior on $\sigma$). Note that if $\mathbf{X}_{(1)}'\mathbf{X}_{\mathbf{l}} = \mathbf{0}$, then $\mathbf{X}_{\mathbf{l}}'\mathbf{Q}_{(1)}\mathbf{X}_{\mathbf{l}} = \mathbf{X}_{\mathbf{l}}'\mathbf{X}_{\mathbf{l}}$, so this would be a standard $g$-type prior.

THEOREM 3. *Under Conditions 1–3 the best predictive model under squared error loss minimizes*

$$(24) \qquad R(M_{\mathbf{l}}) = C + \gamma b^2 \mathbf{w}' \left( \mathbf{P}_{\mathbf{l}} - 2 \sum_{\mathbf{l}^*} p_{\mathbf{l}^*} \mathbf{P}_{\mathbf{l} \cdot \mathbf{l}^*} \right) \mathbf{w},$$

*where* $\mathbf{w} = \mathbf{Q}^{1/2} \mathbf{y}$, $C$ *is a constant and* $\mathbf{l} \cdot \mathbf{l}^*$ *is the dot-product of* $\mathbf{l}$ *and* $\mathbf{l}^*$.

PROOF. Write $\mathbf{x}^* = (\mathbf{x}_{(1)}^*, \mathbf{x}_{(2)}^*)$ and $\mathbf{X} = (\mathbf{X}_{(1)}, \mathbf{X}_{(2)})$, and define $\mathbf{U} = (\mathbf{X}_{(1)}'\mathbf{X}_{(1)})^{-1}\mathbf{X}_{(1)}'$ and $\mathbf{V}_{\mathbf{l}} = (\mathbf{X}_{\mathbf{l}}'\mathbf{Q}_{(1)}\mathbf{X}_{\mathbf{l}})^{-1}\mathbf{X}_{\mathbf{l}}'\mathbf{Q}_{(1)}$. Note that the noncommon variables in $M_{\mathbf{l}}$ are $\mathbf{x}_{(2)}^*\mathbf{H}_{\mathbf{l}_2}$, where $\mathbf{H}_{\mathbf{l}_2}$ is the matrix consisting of the rows of $\mathbf{H}_{\mathbf{l}}$ from $k_1 + 1$ to $k$. With this notation note that

$$(25) \qquad \widehat{y}_{\mathbf{l}}^* = \mathbf{x}_{(1)}^* \widetilde{\boldsymbol{\beta}}_{(1)} + \mathbf{x}_{(2)}^* \mathbf{H}_{\mathbf{l}_2} \widetilde{\boldsymbol{\beta}}_{\mathbf{l}}.$$

Using Condition 3, it is straightforward to show that

$$\mathbf{E}[\boldsymbol{\beta}_{(1)}|\mathbf{y}, \boldsymbol{\beta}_{\mathbf{l}}] = \mathbf{U}(\mathbf{y} - \mathbf{X}_{\mathbf{l}} \boldsymbol{\beta}_{\mathbf{l}}),$$

so that

$$\widetilde{\boldsymbol{\beta}}_{(1)} = \mathbf{E}[\boldsymbol{\beta}_{(1)}|\mathbf{y}] = \mathbf{U}(\mathbf{y} - \mathbf{X}_{\mathbf{l}} \mathbf{E}[\boldsymbol{\beta}_{\mathbf{l}}|\mathbf{y}]).$$

Using this in (25), together with Condition 2, yields

$$\widehat{y}_{\mathbf{l}}^* = \mathbf{x}_{(1)}^* \mathbf{U}(\mathbf{I} - b\mathbf{X}_{\mathbf{l}}\mathbf{V}_{\mathbf{l}})\mathbf{y} + b\mathbf{x}_{(2)}^* \mathbf{H}_{\mathbf{l}_2} \mathbf{V}_{\mathbf{l}} \mathbf{y}$$

$$= \mathbf{x}^* \begin{pmatrix} \mathbf{U}(\mathbf{I} - b\mathbf{X}_{\mathbf{l}}\mathbf{V}_{\mathbf{l}}) \\ b\mathbf{H}_{\mathbf{l}_2} \mathbf{V}_{\mathbf{l}} \end{pmatrix} \mathbf{y}.$$

Defining

$$\mathbf{W}_{\mathbf{l}} = \begin{pmatrix} -\mathbf{U}\mathbf{X}_{\mathbf{l}} \\ \mathbf{H}_{\mathbf{l}_2} \end{pmatrix} \mathbf{V}_{\mathbf{l}},$$

it follows (using Condition 1 in the third equality) that

$$(26) \qquad \begin{aligned} R(M_{\mathbf{l}}) &= \mathbf{E}^{\mathbf{x}^*}[\widehat{y}_{\mathbf{l}}^* - \bar{y}^*]^2 \\ &= b^2 \mathbf{E}^{\mathbf{x}^*} \left[ \mathbf{x}^* \left( \mathbf{W}_{\mathbf{l}} - \sum_{\mathbf{l}^*} p_{\mathbf{l}^*} \mathbf{W}_{\mathbf{l}^*} \right) \mathbf{y} \right]^2 \\ &= \gamma b^2 \mathbf{y}' \left( \mathbf{W}_{\mathbf{l}} - \sum_{\mathbf{l}^*} p_{\mathbf{l}^*} \mathbf{W}_{\mathbf{l}^*} \right)' \mathbf{X}'\mathbf{X} \left( \mathbf{W}_{\mathbf{l}} - \sum_{\mathbf{l}^*} p_{\mathbf{l}^*} \mathbf{W}_{\mathbf{l}^*} \right) \mathbf{y}. \end{aligned}$$

Note that

$$\mathbf{X}\mathbf{W}_{\mathbf{l}} = (\mathbf{X}_{(1)}, \mathbf{X}_{(2)}) \begin{pmatrix} -\mathbf{U}\mathbf{X}_{\mathbf{l}} \\ \mathbf{H}_{\mathbf{l}_2} \end{pmatrix} \mathbf{V}_{\mathbf{l}}$$

$$= -\mathbf{X}_{(1)}\mathbf{U}\mathbf{X}_{\mathbf{l}}\mathbf{V}_{\mathbf{l}} + \mathbf{X}_{\mathbf{l}}\mathbf{V}_{\mathbf{l}} = \mathbf{Q}_{(1)}^{1/2}\mathbf{P}_{\mathbf{l}}\mathbf{Q}_{(1)}^{1/2}.$$

Together with (26) this yields

$$(27) \qquad \begin{aligned} R(M_{\mathbf{l}}) &= \gamma b^2 \mathbf{y}' \mathbf{Q}_{(1)}^{1/2} \left( \mathbf{P}_{\mathbf{l}} - \sum_{\mathbf{l}^*} p_{\mathbf{l}^*} \mathbf{P}_{\mathbf{l}^*} \right) \mathbf{Q}_{(1)} \left( \mathbf{P}_{\mathbf{l}} - \sum_{\mathbf{l}^*} p_{\mathbf{l}^*} \mathbf{P}_{\mathbf{l}^*} \right) \mathbf{y} \\ &= \gamma b^2 \mathbf{w}' \left( \mathbf{P}_{\mathbf{l}} - \sum_{\mathbf{l}^*} p_{\mathbf{l}^*} \mathbf{P}_{\mathbf{l}^*} \right)^2 \mathbf{w}, \end{aligned}$$

the last step utilizing (23) and the fact that $\mathbf{Q}_{(1)}$ is idempotent.

Because $\mathbf{P}_{\mathbf{l}}$ is the projection onto the columns of $\mathbf{Q}_{(1)}^{1/2}\mathbf{X}_{\mathbf{l}}$ that correspond to the nonzero elements of $\mathbf{l}$, it is clear that $\mathbf{P}_{\mathbf{l}}^2 = \mathbf{P}_{\mathbf{l}}$ and $\mathbf{P}_{\mathbf{l}}\mathbf{P}_{\mathbf{l}^*} = \mathbf{P}_{\mathbf{l} \cdot \mathbf{l}^*}$. Expanding the quadratic in (27) with

$$C = \gamma b^2 \mathbf{w}' \left( \sum_{\mathbf{l}^*} p_{\mathbf{l}^*} \mathbf{P}_{\mathbf{l}^*} \right)^2 \mathbf{w}$$

yields the result.  □

COROLLARY 3 (Semi-orthogonal case). *Suppose Conditions* 1–3 *hold and that* $\mathbf{X}_{(2)}'\mathbf{Q}_{(1)}\mathbf{X}_{(2)}$ *is diagonal with positive entries, where the full design matrix is* $\mathbf{X} = (\mathbf{X}_{(1)}\mathbf{X}_{(2)})$. *Then, if the class of models under consideration has graphical model structure, the best predictive model is the median probability model given by* (7).

PROOF. Writing $\mathbf{X}_{(2)}'\mathbf{Q}_{(1)}\mathbf{X}_{(2)} = \mathbf{D}(\mathbf{d})$, the diagonal matrix with diagonal elements $d_i > 0$, $\mathbf{P}_{\mathbf{l}}$ can be expressed as

$$\begin{aligned} \mathbf{P}_{\mathbf{l}} &= \mathbf{Q}_{(1)}^{1/2}\mathbf{X}_{(2)}\mathbf{H}_{\mathbf{l}_2}(\mathbf{H}_{\mathbf{l}_2}'\mathbf{X}_{(2)}'\mathbf{Q}_{(1)}\mathbf{X}_{(2)}\mathbf{H}_{\mathbf{l}_2})^{-1}\mathbf{H}_{\mathbf{l}_2}'\mathbf{X}_{(2)}'\mathbf{Q}_{(1)}^{1/2} \\ &= \mathbf{Q}_{(1)}^{1/2}\mathbf{X}_{(2)}\mathbf{H}_{\mathbf{l}_2}(\mathbf{H}_{\mathbf{l}_2}'\mathbf{D}(\mathbf{d})\mathbf{H}_{\mathbf{l}_2})^{-1}\mathbf{H}_{\mathbf{l}_2}'\mathbf{X}_{(2)}'\mathbf{Q}_{(1)}^{1/2} \\ &= \mathbf{Q}_{(1)}^{1/2}\mathbf{X}_{(2)}(\mathbf{D}(\mathbf{d} \cdot \mathbf{l}))^{-1}\mathbf{X}_{(2)}'\mathbf{Q}_{(1)}^{1/2}. \end{aligned}$$

Hence, defining $\mathbf{u} = \mathbf{X}_{(2)}'\mathbf{Q}_{(1)}^{1/2}\mathbf{w}$, (24) becomes

$$\begin{aligned} R(M_{\mathbf{l}}) &= C + \gamma b^2 \mathbf{u}' \left[ \mathbf{D}(\mathbf{d} \cdot \mathbf{l})^{-1} - 2\sum_{\mathbf{l}^*} p_{\mathbf{l}^*} \mathbf{D}(\mathbf{d} \cdot \mathbf{l} \cdot \mathbf{l}^*)^{-1} \right] \mathbf{u} \\ &= C + \gamma b^2 \sum_{i=1}^{k} u_i^2 d_i^{-1} l_i \left( 1 - 2 \sum_{\mathbf{l}^* : l_i^* = 1} p_{\mathbf{l}^*} \right), \end{aligned}$$

and the conclusion is immediate.  □

Note that $\mathbf{X}'_{(2)}\mathbf{Q}_{(1)}\mathbf{X}_{(2)}$ will be diagonal if either (i) $\mathbf{X}'\mathbf{X}$ is diagonal or (ii) $\mathbf{X}'_{(2)}\mathbf{X}_{(2)}$ is diagonal and $\mathbf{X}'_{(1)}\mathbf{X}_{(2)} = \mathbf{0}$.

COROLLARY 4 (Nested case).  *Suppose Conditions 1–3 hold and that the $M_{\mathbf{l}(j)}, j = 0, \ldots, k$, are a nested sequence of models. Then the best predictive model is the median probability model given by* (7) *or* (9).

PROOF.  For the nested case (24) becomes

$$(28) \quad R(M_{\mathbf{l}(j)}) = C + \gamma b^2 \mathbf{w}' \left( \mathbf{P}_{\mathbf{l}(j)} - 2 \sum_{i=0}^{j-1} p_{\mathbf{l}(i)} \mathbf{P}_{\mathbf{l}(i)} - 2 \sum_{i=j}^{k} p_{\mathbf{l}(i)} \mathbf{P}_{\mathbf{l}(j)} \right) \mathbf{w}.$$

It follows that

$$\begin{aligned}
R(M_{\mathbf{l}(j+1)}) &- R(M_{\mathbf{l}(j)}) \\
&= \gamma b^2 \left( 1 - 2 \sum_{i=j+1}^{k} p_{\mathbf{l}(i)} \right) \mathbf{w}'(\mathbf{P}_{\mathbf{l}(j+1)} - \mathbf{P}_{\mathbf{l}(j)}) \mathbf{w}.
\end{aligned}$$

Since $\mathbf{w}'(\mathbf{P}_{\mathbf{l}(j+1)} - \mathbf{P}_{\mathbf{l}(j)})\mathbf{w} > 0$ and the $(1 - 2\sum_{i=j+1}^{k} p_{\mathbf{l}(i)})$ are increasing in $j$ from $-1$ to $+1$, moving to a larger model will reduce the risk until $(1 - 2\sum_{i=j+1}^{k} p_{\mathbf{l}(i)})$ first turns positive. The conclusion is immediate.  □

EXAMPLE 3.  Consider Hald's regression data [Draper and Smith (1981)], consisting of $n = 13$ observations on a dependent variable $y$ with four potential regressors: $x_1, x_2, x_3, x_4$. Suppose that the following nested models, all including a constant term $c$, are under consideration:

$$\begin{aligned}
M_{\mathbf{l}(1)} &: \{c, x_4\}, & M_{\mathbf{l}(2)} &: \{c, x_1, x_4\}, \\
M_{\mathbf{l}(3)} &: \{c, x_1, x_3, x_4\}, & M_{\mathbf{l}(4)} &: \{c, x_1, x_2, x_3, x_4\},
\end{aligned}$$

again using the notation in (8). We choose the reference prior (14) for the parameters of each model, which effectively means we are using least squares estimates for the predictions and ensures that Conditions 2 and 3 are satisfied. (Here, the models have two common parameters, the constant term and the parameter corresponding to variable $x_4$.) Choosing $\mathbf{Q} = \mathbf{X}'\mathbf{X}$, it follows that the posterior predictive loss of each model is given by (24).

Two choices of model prior probabilities are considered, $P(M_{\mathbf{l}(i)}) = 1/4$, $i = 1, 2, 3, 4$, and $P^*(M_{\mathbf{l}(i)}) = i^{-1}/\sum_{j=1}^{4} j^{-1}$ [the latter type of choice being discussed in, e.g., Jeffreys (1961)]. Default posterior probabilities of each model are then obtained using the Encompassing Arithmetic Intrinsic Bayes Factor, recommended in Berger and Pericchi (1996a, b) for linear models. The resulting model posterior probabilities, $P(M_{l(i)}|\mathbf{y})$ and $P^*(M_{l(i)}|\mathbf{y})$, for the two choices of prior probabilities, respectively, are given in Table 7. The table also presents the normalized posterior predictive loss $R(M_{l(i)}) - C$ for each model.

Since these models are nested, Corollary 4 ensures that the median probability model is the optimal predictive model. Using (9), it is clear from Table 7 that $M_{l(3)}$ is the median probability model for both choices of prior probabilities. And, indeed, the posterior predictive loss of $M_{l(3)}$ is the smallest. Note that $M_{l(3)}$ is the maximum probability model for the first choice of prior probabilities, while $M_{l(2)}$ (which is suboptimal) is the maximum probability model for the second choice.

**5. A geometric formulation.** It was stated in the Introduction that, in general, knowing only the model posterior probabilities does not allow one to determine the optimal predictive model. This is best seen by looking at the problem from a geometric perspective which, furthermore, provides considerable insight into the problem.

Assuming the matrix $\mathbf{Q}$ in (5) is nonsingular and positive definite, consider its Cholesky decomposition $\mathbf{Q} = \mathbf{A}'\mathbf{A}$, where $\mathbf{A}$ is a $k \times k$ upper triangular matrix. The expected posterior loss (16) to be minimized can then be written as

$$(29) \qquad R(M_l) = (\boldsymbol{\alpha}_l - \bar{\boldsymbol{\alpha}})'(\boldsymbol{\alpha}_l - \bar{\boldsymbol{\alpha}}),$$

where $\boldsymbol{\alpha}_l = \mathbf{A}\mathbf{H}_l\tilde{\boldsymbol{\beta}}_l$ is a $k$-dimensional vector and $\bar{\boldsymbol{\alpha}} = \mathbf{A}\bar{\boldsymbol{\beta}} = \sum_l p_l \mathbf{A}\mathbf{H}_l\tilde{\boldsymbol{\beta}}_l$. [If $\mathbf{Q} = \mathbf{X}'\mathbf{X}$ and $\tilde{\boldsymbol{\beta}}_l = \hat{\boldsymbol{\beta}}_l$, one can define $\boldsymbol{\alpha}_l$ as $\boldsymbol{\alpha}_l = \mathbf{X}\mathbf{H}_l\tilde{\boldsymbol{\beta}}_l = \mathbf{X}_l(\mathbf{X}_l'\mathbf{X}_l)^{-1}\mathbf{X}_l'\mathbf{y}$, the projection of $\mathbf{y}$ on the space spanned by the columns of $\mathbf{X}_l$.] It follows that the preferred model will be the one whose corresponding $\boldsymbol{\alpha}_l$ is nearest to $\bar{\boldsymbol{\alpha}}$ in terms of Euclidean distance.

The geometric formulation of the predictive problem follows by representing each model $M_l$ by the point $\boldsymbol{\alpha}_l$. The collection of models thus becomes a collection of points in $k$-dimensional space. The convex hull of these points is a polygon representing the set of possible model averaged estimates $\bar{\boldsymbol{\alpha}}$ as the $p_1$ vary over their range. Hence any point in this polygon is a possible optimal predictive model, depending on the $p_1$, and the goal is to geometrically characterize when each single model is optimal, given that a single model must be used.

TABLE 7
*Posterior probabilities and predictive losses for Hald's data*

|  | $M_{l(1)}$ | $M_{l(2)}$ | $M_{l(3)}$ | $M_{l(4)}$ |
|---|---|---|---|---|
| $P(M_{l(i)}|\mathbf{y})$ | 0.0002 | 0.3396 | 0.5040 | 0.1562 |
| $R(M_{l(i)}) - C$ | 0 | $-808.81$ | $-816.47$ | $-814.43$ |
| $P^*(M_{l(i)}|\mathbf{y})$ | 0.0005 | 0.4504 | 0.4455 | 0.1036 |
| $R(M_{l(i)}) - C$ | 0 | $-808.32$ | $-810.67$ | $-808.31$ |

Consider the simple situation in which we have two covariates $x_1$ and $x_2$ and three possible models:

$$M_{10}: \{x_1\}, \qquad M_{01}: \{x_2\}, \qquad M_{11}: \{x_1, x_2\},$$

again writing, for example, $M_{01}$ instead of $M_{(0,1)}$. These can be represented as three points in the plane. (If the three models had a constant, or intercept, term, then the three points would lie on a plane in three-dimensional space, and the situation would be essentially the same.)

Depending on the sample correlation structure, the triangle whose vertices are $\boldsymbol{\alpha}_{01}$, $\boldsymbol{\alpha}_{10}$ and $\boldsymbol{\alpha}_{11}$ can have three interesting distinct forms. These three forms are plotted in Figure 1. Subregions within each plot will be denoted by the vertices; thus, in Figure 1(a) $[\boldsymbol{\alpha}_{01}, F, C]$ denotes the triangle whose vertices are $\boldsymbol{\alpha}_{01}$, $F$ and $C$.

Each triangle can be divided into optimality subregions, namely the set of those $\bar{\boldsymbol{\alpha}}$ which are closest to one of the $\boldsymbol{\alpha}_1$. These are the regions defined by the solid lines. Thus, in Figure 1(a), the triangle $[\boldsymbol{\alpha}_{10}, F, C]$ defines those points that are closer to $\boldsymbol{\alpha}_{10}$ than to the other two vertices; hence, if $\bar{\boldsymbol{\alpha}}$ were to fall in this region the optimal single model would be $M_{10}$. If $\bar{\boldsymbol{\alpha}}$ were to fall in the triangle $[\boldsymbol{\alpha}_{01}, B, E]$ the optimal single model would be $M_{01}$ and, if $\bar{\boldsymbol{\alpha}}$ were to fall in the region between the two solid lines the optimal single model would be $M_{11}$. It is easy to see that these optimality regions are formed by either (i) connecting the perpendicular bisectors of the sides of the triangle, if all angles are less than or equal to $90°$, or (ii) drawing the perpendicular bisectors of the adjacent side of an angle that is greater than $90°$.

In each plot, $A$, $B$ and $C$ are the midpoints of the line segments $\overrightarrow{\boldsymbol{\alpha}_{10}\boldsymbol{\alpha}_{01}}$, $\overrightarrow{\boldsymbol{\alpha}_{01}\boldsymbol{\alpha}_{11}}$ and $\overrightarrow{\boldsymbol{\alpha}_{10}\boldsymbol{\alpha}_{11}}$, respectively, while $O$ is the midpoint of the triangle. These are of interest because they define regions such that, if $\bar{\boldsymbol{\alpha}}$ lies in the region, then the model corresponding to the vertex in the region has the largest posterior probability. Thus, in Figure 1(a), if $\bar{\boldsymbol{\alpha}}$ lies in the polygon $[\boldsymbol{\alpha}_{10}, A, O, C]$, then $M_{10}$ must be the maximum posterior probability model.

Note that the maximum posterior probability regions do not coincide with the optimal predictive model regions. As a dramatic illustration of the difference, consider Figure 1(a) and suppose that $\bar{\boldsymbol{\alpha}}$ lies on the line segment $\overrightarrow{EF}$. Then $M_{11}$ is the optimal predictive model, even though it has posterior probability 0. Also, either $M_{10}$ or $M_{01}$ has posterior probability at least $1/2$ on this line segment, yet neither is the best predictive model.

The dashed lines form the boundaries defining the median probability models. Thus, if $\bar{\boldsymbol{\alpha}}$ lies in the triangle $[\boldsymbol{\alpha}_{10}, A, C]$, then $M_{10}$ will be the median probability model, while if $\bar{\boldsymbol{\alpha}}$ lies in the polygon $[C, A, B, \boldsymbol{\alpha}_{11}]$, then $M_{11}$ will be the median probability model. To see why this is so, note that the line segment $\overrightarrow{AC}$ consists of the points for which $P(M_{10}|\mathbf{y}) = 1/2$ (i.e., for which $M_{10}$ has posterior probability $1/2$). But then clearly the inclusion probability for variable $x_2$ is also equal to $1/2$ on this line segment, since $p_2 = P(M_{01}|\mathbf{y}) + P(M_{11}|\mathbf{y}) = 1 - P(M_{10}|\mathbf{y})$. Similarly, $\overrightarrow{AB}$ consists of the points for which the inclusion probability for variable $x_1$ is equal to $1/2$. It is immediate that the median probability model in (7) is defined by the indicated regions.

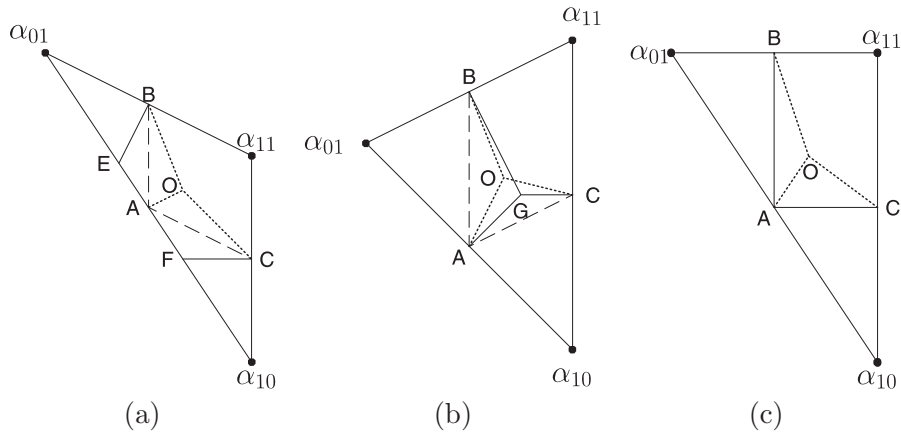

Fig. 1. *Three possible scenarios for the graphical representation of predictive model selection from among $M_{10}$: $\{x_1\}$, $M_{01}$: $\{x_2\}$ and $M_{11}$: $\{x_1, x_2\}$.*

Figures 1(a) and (b) thus show that the median probability model will not always equal the optimal predictive model. Indeed, the two are the same only in the situation of Figure 1(c). In a sense, the theory in the preceding sections arose out of efforts to characterize situations in which the predictive risk representation would be as in Figure 1(c). We found that this is so if $\mathbf{X}'\mathbf{X}$ is diagonal and (17) holds, as in Section 3.1. We also found this to be true in the nested model case discussed in Section 3.3. [Indeed, the resulting figure is simply a rotated version of Figure 1(c).] Subsequently we were able to develop the more general algebraic theories in those sections, but they were based on insights obtained through the geometric formulation.

TABLE 8
*Posterior probabilities and posterior expected losses for Hald's data*

| Model | $P(M_1|\mathbf{y})$ | $R(M_1)$ | Model | $P(M_1|\mathbf{y})$ | $R(M_1)$ |
|---|---|---|---|---|---|
| $c$ | 0.000003 | 2652.44 | $c, 2, 3$ | 0.000229 | 353.72 |
| $c, 1$ | 0.000012 | 1207.04 | $c, 2, 4$ | 0.000018 | 821.15 |
| $c, 2$ | 0.000026 | 854.85 | $c, 3, 4$ | 0.003785 | 118.59 |
| $c, 3$ | 0.000002 | 1864.41 | $c, 1, 2, 3$ | 0.170990 | 1.21 |
| $c, 4$ | 0.000058 | 838.20 | $c, 1, 2, 4$ | 0.190720 | 0.18 |
| $c, 1, 2$ | 0.275484 | 8.19 | $c, 1, 3, 4$ | 0.159959 | 1.71 |
| $c, 1, 3$ | 0.000006 | 1174.14 | $c, 2, 3, 4$ | 0.041323 | 20.42 |
| $c, 1, 4$ | 0.107798 | 29.73 | $c, 1, 2, 3, 4$ | 0.049587 | 0.47 |

One can seek alternative theories based on observations in the geometric formulation. For instance, notice that if the triangle in Figure 1(b) were equilateral, then $O$ and $G$ would coincide and the maximum probability model would equal the optimal predictive model. Unfortunately, we could not find any useful general conditions under which the triangle would be equilateral.

## 6. Concluding comments.

6.1. *When the theory does not apply.* The conditions of the optimality theory for the median probability model are quite strong and will often not apply. Nevertheless, the fact that *only* the median probability model seems to have any optimality theory whatsoever suggests that it might quite generally be successful, even when the optimality theory does not apply.

EXAMPLE 3 (Continued). Suppose that all models (including at least the constant term) are considered for Hald's data. This does not formally satisfy the theory in Section 4, since the models are not nested and the conditions of Theorem 3 do not apply. But here the situation is simple enough that we can directly compute the posterior predictive losses corresponding to each of the possible models, using (16) and assuming equal prior probabilities of the models. The results are given in Table 8.

Computation of the posterior inclusion probabilities yields

$$p_1 = \sum_{1:\, l_1=1} P(M_l|\mathbf{y}) = 0.954556, \qquad p_2 = \sum_{1:\, l_2=1} P(M_l|\mathbf{y}) = 0.728377,$$

$$p_3 = \sum_{1:\, l_3=1} P(M_l|\mathbf{y}) = 0.425881, \qquad p_4 = \sum_{1:\, l_4=1} P(M_l|\mathbf{y}) = 0.553248.$$

Thus the median probability model is $\{c, x_1, x_2, x_4\}$, which from Table 8 clearly coincides with the optimal predictive model. Note that the risk of the maximum probability model $\{c, x_1, x_2\}$ is considerably higher than that of the median probability model.

This example is typical; in our experience the median probability model considerably outperforms the maximum probability model in terms of predictive performance. At the very least this suggests that the median probability model should routinely be determined and reported along with the maximum probability model.

6.2. *When does the median probability model fail?* Suppose that the only models entertained are those with a constant term and a single covariate $x_i$, $i = 1, \ldots, k$, with $k \geq 3$, as well as the model with only a constant term. All models have equal prior probability of $1/(k+1)$. Furthermore, suppose that all covariates are nearly perfectly correlated with each other and with $y$. Then the posterior probability of the constant model will be near zero, and that of each of the other models will coincide with the posterior inclusion probabilities of each of the $x_i$, and will be approximately $1/k$. Since these posterior inclusion probabilities are less than $1/2$, the median probability model will be the constant model, which will have very poor predictive performance compared to any of the other models.

One might be tempted to conclude from this that the median probability model might be problematical if there are highly correlated covariates. We have not yet observed such a difficulty in practice, however. Indeed, the Hald example given in the previous section is an example in which there is high correlation between covariates, yet we saw that the median probability model was still the best.

6.3. *Use of posterior inclusion probabilities.* In addition to being key to defining the median probability model, the posterior inclusion probabilities in (6) can be important tools in assessing the effect of covariates, as indicated in the ANOVA example. One can, furthermore, define joint posterior inclusion probabilities of covariates; these can be very useful in unraveling the effects on model selection of correlations among covariates. See Nadal ([1999](#)) for examples. Finally, posterior inclusion probabilities are a key element in some of the most effective search strategies in model space; compare Berger and Molina ([2002](#)). The importance of posterior inclusion probabilities was emphasized in Mitchell and Beauchamp ([1988](#)).

**Acknowledgments.** The authors are grateful to Merlise Clyde for suggesting the example in Section 6.2. We are also grateful to Edward George for useful discussions.

## REFERENCES

BERGER, J. O. (1997). Bayes factors. In *Encyclopedia of Statistical Sciences, Update* (S. Kotz, C. B. Read and D. L. Banks, eds.) **3** 20–29. Wiley, New York. MR1469744

BERGER, J. O., GHOSH, J. K. and MUKHOPADHYAY, N. (2003). Approximations and consistency of Bayes factors as model dimension grows. *J. Statist. Plann. Inference* **112** 241–258. MR1961733

BERGER, J. O. and MOLINA, M. (2002). Discussion of "A case study in model selection," by K. Viele, R. Kass, M. Tarr, M. Behrmann and I. Gauthier. In *Case Studies in Bayesian Statistics VI. Lecture Notes in Statist.* **167** 112–125. Springer, New York. MR1955885

BERGER, J. O. and PERICCHI, L. R. (1996a). The intrinsic Bayes factor for model selection and prediction. *J. Amer. Statist. Assoc.* **91** 109–122. MR1394065

BERGER, J. O. and PERICCHI, L. R. (1996b). The intrinsic Bayes factor for linear models. In *Bayesian Statistics 5* (J. M. Bernardo, J. O. Berger, A. P. Dawid and A. F. M. Smith, eds.) 25–44. Oxford Univ. Press. MR1425398

BERGER, J. O. and PERICCHI, L. R. (2001). Objective Bayesian methods for model selection: Introduction and comparison (with discussion). In *Model Selection* (P. Lahiri, ed.) 135–207. IMS, Beachwood, OH.

BERGER, J. O., PERICCHI, L. R. and VARSHAVSKY, J. (1998). Bayes factors and marginal distributions in invariant situations. *Sankhyā Ser. A* **60** 307–321. MR2000753

BERNARDO, J. and SMITH, A. F. M. (1994). *Bayesian Theory*. Wiley, New York. MR1274699

CHIPMAN, H., GEORGE, E. I. and McCULLOCH, R. E. (2001). The practical implementation of Bayesian model selection (with discussion). In *Model Selection* (P. Lahiri, ed.) 66–134. IMS, Beachwood, OH. MR2000752

CLYDE, M. A. (1999). Bayesian model averaging and model search strategies. In *Bayesian Statistics 6* (J. M. Bernardo, J. O. Berger, A. P. Dawid and A. F. M. Smith, eds.) 157–185. Oxford Univ. Press. MR1723497

CLYDE, M. A., DESIMONE, H. and PARMIGIANI, G. (1996). Prediction via orthogonalized model mixing. *J. Amer. Statist. Assoc.* **91** 1197–1208.

CLYDE, M. A. and GEORGE, E. I. (1999). Empirical Bayes estimation in wavelet nonparametric regression. In *Bayesian Inference in Wavelet-Based Models. Lecture Notes in Statist.* **141** 309–322. Springer, New York. MR1699849

CLYDE, M. A. and GEORGE, E. I. (2000). Flexible empirical Bayes estimation for wavelets. *J. R. Stat. Soc. Ser. B Stat. Methodol.* **62** 681–698. MR1796285

CLYDE, M. A. and PARMIGIANI, G. (1996). Orthogonalizations and prior distributions for orthogonalized model mixing. In *Modelling and Prediction* (J. C. Lee, W. Johnson and A. Zellner, eds.) 206–227. Springer, New York.

CLYDE, M. A., PARMIGIANI, G. and VIDAKOVIC, B. (1998). Multiple shrinkage and subset selection in wavelets. *Biometrika* **85** 391–401. MR1431069

DRAPER, N. and SMITH, H. (1981). *Applied Regression Analysis*, 2nd ed. Wiley, New York. MR610978

GEORGE, E. I. and FOSTER, D. P. (2000). Calibration and empirical Bayes variable selection. *Biometrika* **87** 731–747. MR1813972

GEORGE, E. I. and McCULLOCH, R. E. (1993). Variable selection via Gibbs sampling. *J. Amer. Statist. Assoc.* **88** 881–889.

JEFFREYS, H. (1961). *Theory of Probability*, 3rd ed. Oxford Univ. Press.

HOETING, J. A., MADIGAN, D., RAFTERY, A. E. and VOLINSKY, C. T. (1999). Bayesian model averaging: A tutorial (with discussion). *Statist. Sci.* **14** 382–417. MR187257

MITCHELL, T. J. and BEAUCHAMP, J. J. (1988). Bayesian variable selection in linear regression (with discussion). *J. Amer. Statist. Assoc.* **83** 1023–1036. MR997578

MONTGOMERY, D. C. (1991). *Design and Analysis of Experiments*, 3rd ed. Wiley, New York. MR1076621

MUKHOPADHYAY, N. (2000). Bayesian and empirical Bayesian model selection. Ph.D. dissertation, Purdue Univ.

MÜLLER, P. (1999). Simulation-based optimal design. In *Bayesian Statistics 6* (J. M. Bernardo, J. O. Berger, A. P. Dawid and A. F. M. Smith, eds.) 459–474. Oxford Univ. Press. MR1723509

NADAL, N. (1999). El análisis de varianza basado en los factores de Bayes intrínsecos. Ph.D. thesis, Univ. Simón Bolívar, Venezuela.

O'Hagan, A. (1995). Fractional Bayes factors for model comparison (with discussion). *J. Roy. Statist. Soc. Ser. B* **57** 99–138. MR1325379

Scheffé, H. (1959). *The Analysis of Variance.* Wiley, New York. MR116429

Shibata, R. (1983). Asymptotic mean efficiency of a selection of regression variables. *Ann. Inst. Statist. Math.* **35** 415–423. MR739383

Zellner, A. (1986). On assessing prior distributions and Bayesian regression analysis with *g*-prior distributions. In *Bayesian Inference and Decision Techniques*: *Essays in Honor of Bruno de Finetti* (P. K. Goel and A. Zellner, eds.) 233–243. North-Holland, Amsterdam. MR881437

Dipartimento di Economia
Università Roma Tre
via Ostiense 139
00154 Roma
Italy
e-mail: marilena.barbieri@uniroma3.it

Institute of Statistics
  and Decision Sciences
Duke University
Durham, North Carolina 27708-0251
USA
e-mail: berger@stat.duke.edu
url: www.stat.duke.edu/~berger